\documentclass[12pt]{amsart}
\usepackage{amscd,amssymb,hyperref}
\usepackage{pstricks}

\theoremstyle{plain}

\theoremstyle{remark}

\theoremstyle{definition}

\numberwithin{equation}{section}

\begin{document}
\baselineskip=26pt
 \title {Closed string operators in topology leading to Lie bialgebras and higher string algebra}
\author{Moira Chas and Dennis Sullivan}
\maketitle

Imagine a collection of closed oriented curves depending on
parameters in a smooth d-manifold $M$.  Along a certain locus of
configurations strands of the curve may intersect at certain sites
in $M$. At these sites in $M$ the curves may be cut and
reconnected in some way.  One obtains operators on the set of
parametrized collections of closed curves in $M$.  By making the
coincidences transversal and compactifying, the operators can be
made to act in the algebraic topology of the free loop space of
$M$ when $M$ is oriented.  The process reveals collapsing sub
graph combinatorics like that for removing infinities from Feynman
graphs.

Let $\mathbb{H}$ denote the equivariant homology with rational
coefficients of the pair (Maps ($S^1,M$), constant maps ($S^1,M$))
relative to the $S^1$ action of rotating the source.  $\mathbb{H}$
is called the reduced equivariant homology of the free loop space
of $M$.

Associated to the diagrams I and II of the figure there are
operators $c_n:\mathbb{H}^{\otimes n}\rightarrow \mathbb{H}$ and
$s_n:\mathbb{H} \rightarrow \mathbb{H}^{\otimes n}$, $n=2,3,4,..$
These operators satisfy various relations e.g. $c_2$ and $s_2$
satisfy the relations-Jacobi, coJacobi, and Drinfeld compatibility
which utilize diagrams III, IV, V, and VI. \newline Diagram I (the
operator $c_n$, $n=5$):
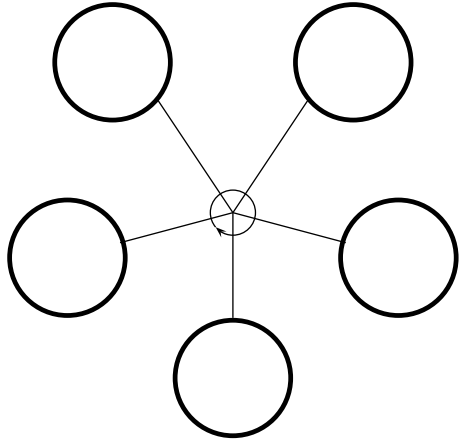
\[
\begin{pspicture}(0,0)(6,6)
% \psgrid(0,0)(6,6)
 \pscircle[linewidth=1.8pt](1.4,5){0.8}
 \pscircle[linewidth=1.8pt](4.6,5){0.8}
 \pscircle[linewidth=1.8pt](3,0.8){0.8}
 \pscircle[linewidth=1.8pt](0.8,2.4){0.8}
 \pscircle[linewidth=1.8pt](5.2,2.4){0.8}
 \psarc[linewidth=0.5pt, arrowsize=3pt]{<-}(3,3){0.3}{220}{220}
 \psline[linewidth=.5pt](3,3)(4,4.5)
 \psline[linewidth=.5pt](3,3)(2,4.5)
 \psline[linewidth=.5pt](3,3)(3,1.6)
 \psline[linewidth=.5pt](3,3)(1.5,2.6)
 \psline[linewidth=.5pt](3,3)(4.5,2.6)
\end{pspicture}
\]

 Diagram II (the operator $s_n$, $n=6$):
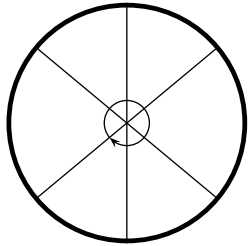
\[
\begin{pspicture}(0,0)(4,4)
% \psgrid(0,0)(4,4)
 \pscircle[linewidth=1.8pt](2,2){1.6}
 \psarc[linewidth=0.5pt, arrowsize=3pt]{<-}(2,2){0.3}{220}{220}
 \psline[linewidth=.5pt](0.8,1)(3.2,3)
 \psline[linewidth=.5pt](3.2,1)(0.8,3)
 \psline[linewidth=.5pt](2,3.6)(2,0.4)
\end{pspicture}
\]
Diagram III (Jacobi identity for $c_2$):
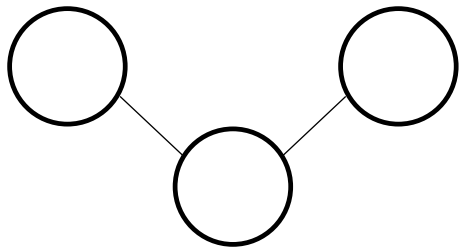
\[
\begin{pspicture}(0,0)(6,3.2)
% \psgrid(0,0)(6,3.2)
 \pscircle[linewidth=1.8pt](3,0.8){0.8}
 \pscircle[linewidth=1.8pt](0.8,2.4){0.8}
 \pscircle[linewidth=1.8pt](5.2,2.4){0.8}
 \psline[linewidth=.5pt](2.35,1.2)(1.5,2)
 \psline[linewidth=.5pt](3.65,1.2)(4.5,2)
\end{pspicture}
\]
 Diagram IV (Cojacobi identity for $s_2$):
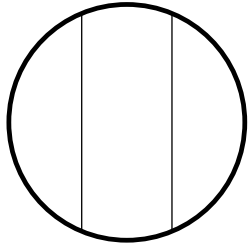
\[
\begin{pspicture}(0,0)(4,4)
% \psgrid(0,0)(4,4)
 \pscircle[linewidth=1.8pt](2,2){1.6}
 \psline[linewidth=.5pt](2.6,.55)(2.6,3.45)
 \psline[linewidth=.5pt](1.4,.55)(1.4,3.45)
\end{pspicture}
\]
 Diagram V (Drinfeld compatibility for $c_2$ and $s_2$):
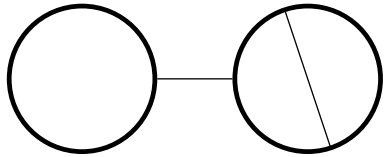
\[
\begin{pspicture}(0,0)(5,3)
% \psgrid(0,0)(5,3)
 \pscircle[linewidth=1.8pt](1,1.5){1}
 \pscircle[linewidth=1.8pt](4,1.5){1}
 \psline[linewidth=.5pt](2,1.5)(3,1.5)
 \psline[linewidth=.5pt](3.7,2.4)(4.3,0.6)
\end{pspicture}
\]
 Diagram VI (Drinfeld compatibility for $c_2$ and $s_2$):
\[
\begin{pspicture}(0,0)(5,3)
% \psgrid(0,0)(5,3)
 \pscircle[linewidth=1.8pt](1,1.5){1}
 \pscircle[linewidth=1.8pt](4,1.5){1}
 \psline[linewidth=.5pt](1.8,1)(3.2,1)
 \psline[linewidth=.5pt](1.8,2)(3.2,2)
\end{pspicture}
\]
Thus we prove the

\textit{Theorem}: The reduced equivariant homology $\mathbb{H}$ of
the free loop space of a smooth oriented d-manifold $M$ has the
structure of a Lie bialgebra generated by string operators
$c_2:\mathbb{H}\otimes \mathbb{H}\rightarrow\mathbb{H}$ and
$s_2:\mathbb{H}\rightarrow\mathbb{H}\otimes\mathbb{H}$ of degree
$2-d$. The string operators $c_n:\mathbb{H}^{\otimes
n}\rightarrow\mathbb{H}$  and $s_n:\mathbb{H}\rightarrow
\mathbb{H}^{\otimes n}$ $S_n$ are also defined for $n>2$ and have
degree $n+ (1-n) d$. The relations conjecturally satisfied by
$c_n$ and $s_n$ are described below.

\textit{Problem and Conjecture}: There is evidence (see $c_n$ and
$s_n$ Identities paragraph below) that one may prove the $c_n$ and
$s_n$ generate an algebraic structure on $\mathbb{H}$ that is
Koszul dual in the sense of [Gan] to the positive boundary version
(the number of inputs and the number of outputs are both positive)
of the algebraic structure (genus zero) of symplectic topology
[M].

 \textit{Remark}: This Lie bialgebra is non trivial. For example,
 when $d=2$ it is isomorphic to the Lie bialgebra discovered by
 Goldman and Turaev[T].

 \textit{General Program}: The theorem is part
of a more elaborate structure of closed string operators including
higher genus acting at the chain level of the spaces of maps of
families of closed curves into $M$. We describe the elements of
this general theory and illustrate with the examples indicated by
diagrams I, II, III, IV, V, and VI which are required to treat the
theorem and diagram VII used to prove the genus one involutive
identity.

\bigskip

\begin{pspicture}(14,5)%\psgrid[gridcolor=gray](14,5)

 \psline[arrowsize=6pt]{->}(6.9,4.7)(7.1,4.7)
\rput(3,4.6){Diagram VII}
 \pscircle(7,2.5){2.2}
\psline(5.4,1)(8.6,4)\psline(5.4,4)(6.8,2.7)
\psline(7.2,2.3)(8.6,1) %(7.2,2.4)
\psdot(6.4,3.08) \rput(4,0) {(Involutive identity $c_2 \cdot
s_2=0$)}

\end{pspicture}

\bigskip

\textit{Generalized Chord Diagram}: First of all there are
operators of order zero, order one, etc.  The operators of order
zero are associated to diagrams as above. The most general
\textit{order zero diagram }$D$ (generalized chord diagram)is
specified by the data- a finite subset $F$ of a union $C$ of
directed circles, a partition of the subset $F$ into parts of
cardinality at least two, and a cyclic order on each part. Higher
order operations correspond to general chord diagrams with
additional combinatorial data (see paragraphs Diagrams V, Diagram
VI). To the order zero diagram one can associate a graph $\Gamma
^\prime(D)$ by attaching an $n$-prong respecting the cyclic order
(as in diagram I or II) to each part of the partition. There is
also a surface $\Sigma^\prime(D)$ obtained by crossing the
collection of circles $C$ with the unit interval [0,1] and
attaching along n segments a thickened (n-prong), a disk with $n$
tabs, to $C$x$\{1\}$. So $\Sigma^\prime (D)$ equals $C$x$[0,1]$
union (thickened n-prongs) along $C$x$\{1\}$.

\textit{Ribbon surface and cyclic graph of a diagram}: We denote
by $\Sigma(D)$ the surface obtained from $\Sigma^\prime (D)$ by
collapsing to points the n-prongs slightly extended into the
collar to level $C$x$\{1/2\}\ $.  The quotient of $C$x$\{1/2\}\ $
by this collapse is called the cyclically ordered graph
$\Gamma(D)$ of the diagram $D$. $\Sigma(D)$ is a \textit{ribbon
surface} associated to the \textit{cyclically ordered graph}
$\Gamma (D)$, associated to the \textit{generalized chord diagram}
$D$.

\textit{Closed String Operator}: Let us try to define the chain
operator associated to $D$.  The input is a bundle $\eta$ with
fibre $C$ and a map of the total space of $\eta$ into $M$ smooth
along the $C$ direction. (These objects can be used to compute
tensor products of the equivariant homology reduced or not).  We
form the associated bundle to $\eta$ with fibre all isotopic
configurations of $F$ embedded in $C$.  We restrict attention to
the \textit{$D$-locus} in the total space of the associated
bundle, defined as the locus where the map identifies in $M$ all
the points in each part of the partition of $F$ (now considered as
a subset of $C$x$\{1/2\})$. Along this locus there are induced
maps of the quotient graph $\Gamma(D)$ of $C$ into $M$ and thus
also of the ribbon surface $\Sigma(D)$ which retracts to
$\Gamma(D)$. The restriction of this map on $\Sigma(D)$ to
$C^\prime$, the rest of the boundary of $\Sigma (D)$ besides $C$
(actually $C$x$\{0\})$ in each fibre is a map of the total space
of a $C^\prime$ fibration over the $D$ locus into $M$. This is the
output of the $D$ operation at the level of set theory or
topological spaces.

\textit{Chain Operator}: To have a good object in algebraic
topology there are several issues compactness, transversality, and
orientation. We want the output to be a relative chain
representing an element in a chain complex computing
$\mathbb{H}^{\otimes i}$, $i=$ number of components of $C$, if the
input is.  We now discuss these issues in the order mentioned.

\textit{Diagram I}: Now the compactness property of chains is not
a problem for diagram I because its configuration space is already
compact being an $n$-torus.  Similarly, compactness is easy to
arrange for any cactus diagram (diagram $III^\prime$) generalizing
diagram III. (see paragraph Diagram III below).  The homological
content of this part of theory was discussed in [CS] where the non
reduced equivariant homology, the Lie bracket $c_2$ and the higher
analogues $c_n$ were considered.

\textit{Diagram II}: For diagram II the configuration space is non
compact and there is an obstacle to overcome.  Note however that
when two prongs come close enough together (relative to the input
C-family) the output $C^\prime$ family has for each parameter at
least one component which is small in $M$.  The role of the
reduced equivariant homology is to take advantage of this fact.
Because families with constant (or small) map components are
considered to be null or zero for the reduced discussion we obtain
that the configuration space of diagram II is (relatively) compact
for the purpose of producing a (relative) chain for the
computation of $\mathbb{H}^{\otimes i}$, the reduced equivariant
homology.

\textit{Relative Compactness}: To treat the chain operator for a
general diagram we need to complete its configuration space enough
so that this relative compactness is achieved.  Namely, any strata
omitted from a true compactification must correspond to output
families which for each parameter have a positive number of tiny
components.

\textit{Constraint Normal Bundle}: Besides this compactness
consideration the other "sine qua non" issue is the normal bundle
to the constraint locus.  Imposing the conditions defining the
locus of the diagram in the cases above amounts to taking the
preimage of a diagonal with a normal bundle.  Then the transversal
preimage of a chain will be a chain (the Thom map at the chain
level) and we can work in the context of algebraic topology.

We refer to this as the constraint normal bundle issue and we must
keep this \textit{constraint normal bundle} as we add pieces to
the configuration spaces to obtain (relative) compactness.

\textit{Diagram III}: Let us consider diagram III. When two
different chords on the middle circle coalesce we converge to
constraints as in diagram I (n=3).  The number of constraints is
the same (two) and the normal bundles fit together perfectly. This
works for the normal bundle consideration to compactify all the
configuration spaces of cactus like diagrams, diagram $III^\prime$
(planar trees with circles inserted at some of the vertices)
leading to genus zero $n$ to one operations.  The point is that
compactness is achieved by allowing \textit{different} parts of
the partitions to coalesce where the constraint normal bundle
persists. Because of the tree like property it is impossible to
let parts self collide. Except for orientations the discussion is
conceptually complete in these cases- we have actual compactness
and a normal bundle to define a transversal Thom chain map, in the
equivariant context. (compare [CS] and [V]).
\[
\begin{pspicture}(14,12)%\psgrid[gridcolor=gray](14,12)

\rput(2,11){Diagram $III^\prime:$}
\pscircle[linewidth=1.5pt](0.6,3.5){0.5}
\pscircle[linewidth=1.5pt](2.5,5){0.7}
\pscircle[linewidth=1.5pt](1.7,6.8){0.5}
\pscircle[linewidth=1.5pt](5,4){1}
\pscircle[linewidth=1.5pt](6.3,1.8){0.5}
\pscircle[linewidth=1.5pt](7.6,6){0.8}
\pscircle[linewidth=1.5pt](7,8.8){0.6}
\pscircle[linewidth=1.5pt](11,6.4){0.6}
\pscircle[linewidth=1.5pt](9.4,8.6){0.5}

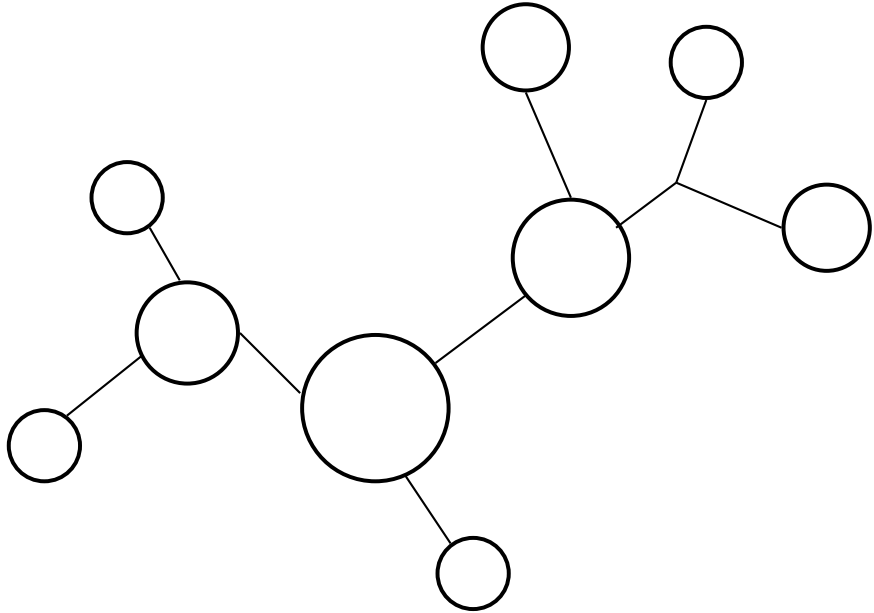
\psline(6,2.2)(5.4,3.1) \psline(4,4.2)(3.2,5)
\psline(0.9,3.9)(1.9,4.7) \psline(2,6.4)(2.4,5.7)
\psline(7,5.5)(5.8,4.6) \psline(7.6,6.8)(7,8.2)
\psline(8.2,6.4)(9,7) \psline(9,7)(9.4,8.1) \psline(9,7)(10.4,6.4)

\end{pspicture}
\]
\[
\begin{pspicture}(14,12)%\psgrid[gridcolor=gray](14,12)
\psellipse(6,5.4)(7,5.4) \psline[linestyle=solid](12,8.2)(9,7)
\psline[linestyle=solid](7.5,10.6)(9,7)
\psline[linestyle=solid](11.2,2)(9,7)
\pscurve[linestyle=solid](7,10.7)(7.6,7.4)(5,10.7)
\pscurve[linestyle=solid](3,10.3)(6.3,5)(10,1)
\pscurve[linestyle=solid](2.4,10)(4,5)(1,1.6)
\pscurve[linestyle=solid](9,0.6)(5.6,2.8)(5,0.1)
\pscurve[linestyle=solid](2,9.8)(2.2,6)(-1,6)
\pscurve[linestyle=solid](-0.75,4)(1.5,4.4)(0.6,2)
\rput(4,11.3){Diagram $IV^\prime:$}
\end{pspicture}
\]
\[
\begin{pspicture}(14,12)%\psgrid[gridcolor=gray](14,12)

 \psellipse(6,5.4)(7,5.4)
 \pscircle[linewidth=1.5pt](0.6,3.5){0.5}
\pscircle[linewidth=1.5pt](2.5,5){0.7}
\pscircle[linewidth=1.5pt](1.7,6.8){0.5}
\pscircle[linewidth=1.5pt](5,4){1}
\pscircle[linewidth=1.5pt](6.3,1.8){0.5}
\pscircle[linewidth=1.5pt](7.6,6){0.8}
\pscircle[linewidth=1.5pt](7,8.8){0.6}
\pscircle[linewidth=1.5pt](11,6.4){0.6}
\pscircle[linewidth=1.5pt](9.4,8.6){0.5} \psline(6,2.2)(5.4,3.1)
\psline(4,4.2)(3.2,5) \psline(0.9,3.9)(1.9,4.7)
\psline(2,6.4)(2.4,5.7) \psline(7,5.5)(5.8,4.6)
\psline(7.6,6.8)(7,8.2) \psline(8.2,6.4)(9,7)
\psline(9,7)(9.4,8.1) \psline(9,7)(10.4,6.4) \rput(4,11.3){Duality
Diagram:} \psline[linestyle=dashed](12,8.2)(9,7)
\psline[linestyle=dashed](7.5,10.6)(9,7)
\psline[linestyle=dashed](11.2,2)(9,7)
\pscurve[linestyle=dashed](7,10.7)(7.6,7.4)(5,10.7)
\pscurve[linestyle=dashed](3,10.3)(6.3,5)(10,1)
\pscurve[linestyle=dashed](2.4,10)(4,5)(1,1.6)
\pscurve[linestyle=dashed](9,0.6)(5.6,2.8)(5,0.1)
\pscurve[linestyle=dashed](2,9.8)(2.2,6)(-1,6)
\pscurve[linestyle=dashed](-0.75,4)(1.5,4.4)(0.6,2)

\end{pspicture}
\]

\textit{Diagram IV}: Now consider diagram IV.  When different
chords come together \textit{at one point} the constraints
converge to those of a diagram like II (n=3).  The normal bundles
fit together perfectly as before. By adding the configuration
spaces of diagram II to that of IV we obtain a \textit{good}
constraint normal bundle and \textit{relative} compactness for
Diagram IV. This also applies to more general internal chord or
n-prong diagrams on one circle (disjoint n-prongs in one circle
Diagram $IV^\prime$) which are dual to planar cactii with roles of
input and output interchanged (see Duality Diagram). We can add
strata preserving the constraint normal bundle and achieve
relative compactness. Thus dual to all the $n$ to one genus zero
operations Diagram $III^\prime$ from the Diagram III paragraph we
have one to $n$ genus zero operations Diagram $IV^\prime$ in the
reduced theory (however see next paragraph).

\textit{Null chains and Degenerate chains}: There is one
additional caveat about working in the reduced theory even for the
diagrams considered up to now.  We need to know the null
subcomplex consisting of maps where at least one component is
constant is invariant by the operations.  At first glance this
seems problematic but it works out in the end. A constant loop
component may be cut up by an operation and mixed into other
components.  If so we no longer have a constant component.
However, the situation is saved because we obtain a degenerate
chain-one whose geometric dimension is too low.  For as the
parameters of the configurations vary in the component which is
mapped to a constant loop the image chain is not varying.  A
transversal pull back will not have the full homological dimension
and can be ignored.  In fact we can mod out by degenerate chains
from the beginning.  (This point has to be considered carefully
when extraordinary homology theories are studied here).  On the
other hand if the operation doesn't touch the null component we
still have a null component in the output.  In summary definable
string chain operators act in the relative complex defining the
reduced equivariant theory. For example we have all the
ingredients now to define the operations $c_n$ and $s_n$,
$n=2,3,4,...$ in the reduced equivariant theory.

One need only add that diagrams I and II correspond to
\textit{cycles} since no strata were added to create (relative)
compactness.  Thus the chain operators corresponding to $c_n$ and
$s_n$ commute with the $\partial$ operators on chains and pass to
homology. (Orientations will be discussed in the paragraph below).

To prove the relations of a Lie bialgebra among the compositions
of $c_2$ and $s_2$ we have to consider diagrams V and VI which
bring forth two further considerations.

\textit{Diagram V}: To achieve relative compactness for the
configurations space of diagram V we have to let the two chords
touch at one point which is a case already considered above for
Diagram III (and IV).  We also have to allow the internal chord of
Diagram V to collapse to the endpoint of the connecting chord
(from opposite sides only- because a one sided approach leads to a
tiny output circle and a null chain).  This creates in the limit a
diagram of type I for $n=2$ and if we do nothing else the number
of constraints goes down and we lose the normal bundle property.
However, the collapsing internal chord and the constraint that
values at the endpoints of this chord coincide say that in the
limit the derivative of the map in the $C$ direction is null at
the limit point.

Thus we are led to an order one diagram, a diagram of order zero
of type I for $n=2$ with the additional data that one of the
attaching points of the chord is a point of multiplicity two. This
means that when the locus of this diagram of order one is defined
the condition coincidence of values at the endpoints of the chord
is augmented by the condition that the 1st derivative in the $C$
direction is zero at the point of multiplicity two.  (In general
at a point of multiplicity $k$ the first ($k-1$) derivatives would
be required to be zero).  Then the constraint normal bundle
extends continuously over the added stratum for the relative
compactness. This treats diagram V.

\textit{Diagram VI}: One more feature appears in treating diagram
VI. To the generic $4D$ configuration space of diagram VI we add
three $3D$ strata and two $2D$ strata for relative compactness.
One of the $3D$ strata will involve a new kind of consideration
similar to that in the Fulton MacPherson compactification of
configuration spaces [FM]. The other strata will be of the type
already considered. Namely, two of the $3D$ strata allow the two
chords to touch on one circle or the other.  This has already been
considered. The new case appears when the two chords approach each
other (on opposite sides again) at both endpoints at commensurable
distances.  We add a $3D$ stratum to our space which records the
limiting single chord (two parameters) and a third parameter which
can record the signed ratio of the small distances in the
approach.  We call this an FM-stratum.

We also add two strata for the chords touching first at one
endpoint and then at the second endpoint which-like the discussion
of diagram V- produces in each case a multiple point. These strata
account for the approach of chords at both endpoints with
incommensurable distances at the endpoints.

When defining the locus for this completion (to a relatively
compact configuration space) we treat the FM stratum in the
following way.  We have a chord diagram of type I ($n=2$) with a
third ratio parameter $\lambda$.  We ask first that the map agrees
at endpoints of the chord as before and then ask further that
derivatives in the $C$ direction at these points be proportional
with ratio $1/\lambda$, (the factor $1/\lambda$ because distances
appear in the denominator when calculating derivatives).

If the other strata are treated as described above the constraint
normal bundle extends continuously over this entire (relative)
compactification of the configuration space of diagram VI.  This
treats diagram VI.

Now we have all the ingredients to define the chain operators of
diagrams I through VI up to a question of orientation.

\textit{Orientations:} It is possible to avoid a nightmare of sign
difficulties using a categorial approach to orientations motivated
by [D]. First there is the "graded line" functor from finite
dimensional real vector spaces to $Z/2$ graded vector spaces.  It
assigns to $V$ the top exterior power placed in even degree if
dimension $V$ is even and in odd degree if dimension $V$ is odd.

An orientation of $V$ is by definition a generator of the graded
line of $V$ up to positive multiples.  For any finite family of
spaces $V_f$, $f \epsilon F$ there is a canonical notion of direct
sum $\underset {f \epsilon F}{\bigoplus}$ $V_f$ and if each $V_f$
is oriented a canonical orientation of $\underset {f \epsilon
F}{\bigoplus}$ $V_f$. Similarly in an exact sequence of spaces $0
\rightarrow V \rightarrow W \rightarrow V^\prime \rightarrow 0$ an
orientation of any two determines canonically an orientation of
the third.

We apply this to the previous constructions as follows.  Fix an
orientation of $M$.  Then any product $\underset{f \epsilon
F}{\Pi} M$ indexed by a finite set has a canonical orientation by
the first fact above.  Using the second fact as well, any diagonal
corresponding to a part of $F$ and its normal bundle each has a
canonical orientation.  Further, if the base of a $C$ bundle is
oriented the total space $\eta$ is canonically oriented (since $C$
is oriented) and the total space of the associated bundle of $F$
configurations in $C$ is oriented (if $F$ is ordered up to even
permutations).  In our examples Diagrams I, II, III, IV, V, VI
this ordering on $F$ comes from the ordering of input or output
components. In example Diagram VII it comes from the ordering of
the chords. Combining all this the transversal pull back is also
canonically oriented. The actual definition of the operators $c_n$
and $s_n$ uses these canonical orientations.

\textit{Lie bialgebra identities, Jacobi}:

The proof of Jacobi for $c_2$ uses several versions of diagram
III, to calculate [[1,2],3]=$c_2$ $(c_2(1 \otimes 2)\otimes 3)$
and its cyclic permutations. The arrow on the chords is determined
by the order of input arguments and determines, via an ordering in
$F$ up to even permutations, the orientation. The cancellation is
indicated.

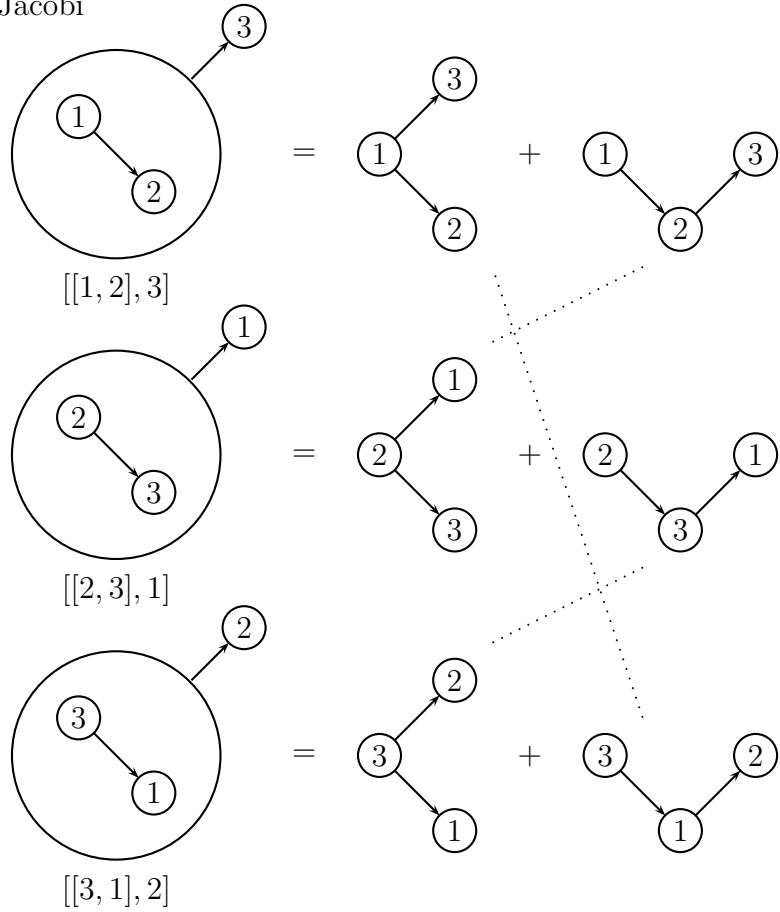
\begin{pspicture}(0,0)(11,12)
 %\psgrid(0,0)(11,12)
 \pscircle(2,2){1.4}
 \pscircle(1.5,2.5){0.3} \rput(1.5,2.5){3}
 \pscircle(2.5,1.5){0.3} \rput(2.5,1.5){1}
 \pscircle(3.7,3.7){0.3} \rput(3.7,3.7){2}
 \psline{->}(1.7,2.3)(2.3,1.7)  \psline{->}(3,3)(3.5,3.5)
\rput(1,12){Jacobi}
 \rput(2,0.2){$[[3,1],2]$}
 \rput(4.5,2){=}  \rput(7.5,2){+}
 \pscircle(5.5,2){0.3} \rput(5.5,2){3}
 \pscircle(6.5,3){0.3} \rput(6.5,3){2}
 \pscircle(6.5,1){0.3} \rput(6.5,1){1}
 \pscircle(8.5,2){0.3} \rput(8.5,2){3}
 \pscircle(9.5,1){0.3} \rput(9.5,1){1}
 \pscircle(10.5,2){0.3} \rput(10.5,2){2}
 \psline{->}(5.7,2.2)(6.3,2.8) \psline{->}(5.7,1.8)(6.3,1.2)
 \psline{->}(8.7,1.8)(9.3,1.2) \psline{->}(9.7,1.2)(10.3,1.8)
 \pscircle(2,6){1.4}
 \pscircle(1.5,6.5){0.3} \rput(1.5,6.5){2}
 \pscircle(2.5,5.5){0.3} \rput(2.5,5.5){3}
 \pscircle(3.7,7.7){0.3} \rput(3.7,7.7){1}
 \psline{->}(1.7,6.3)(2.3,5.7)  \psline{->}(3,7)(3.5,7.5)
 \rput(2,4.2){$[[2,3],1]$}
 \rput(4.5,6){=}  \rput(7.5,6){+}
 \pscircle(5.5,6){0.3} \rput(5.5,6){2}
 \pscircle(6.5,7){0.3} \rput(6.5,7){1}
 \pscircle(6.5,5){0.3} \rput(6.5,5){3}
 \pscircle(8.5,6){0.3} \rput(8.5,6){2}
 \pscircle(9.5,5){0.3} \rput(9.5,5){3}
 \pscircle(10.5,6){0.3} \rput(10.5,6){1}
 \psline{->}(5.7,6.2)(6.3,6.8) \psline{->}(5.7,5.8)(6.3,5.2)
 \psline{->}(8.7,5.8)(9.3,5.2) \psline{->}(9.7,5.2)(10.3,5.8)
 \pscircle(2,10){1.4}
 \pscircle(1.5,10.5){0.3} \rput(1.5,10.5){1}
 \pscircle(2.5,9.5){0.3} \rput(2.5,9.5){2}
 \pscircle(3.7,11.7){0.3} \rput(3.7,11.7){3}
 \psline{->}(1.7,10.3)(2.3,9.7)  \psline{->}(3,11)(3.5,11.5)
 \rput(2,8.2){$[[1,2],3]$}
 \rput(4.5,10){=}  \rput(7.5,10){+}
 \pscircle(5.5,10){0.3} \rput(5.5,10){1}
 \pscircle(6.5,11){0.3} \rput(6.5,11){3}
 \pscircle(6.5,9){0.3} \rput(6.5,9){2}
 \pscircle(8.5,10){0.3} \rput(8.5,10){1}
 \pscircle(9.5,9){0.3} \rput(9.5,9){2}
 \pscircle(10.5,10){0.3} \rput(10.5,10){3}
 \psline{->}(5.7,10.2)(6.3,10.8) \psline{->}(5.7,9.8)(6.3,9.2)
 \psline{->}(8.7,9.8)(9.3,9.2) \psline{->}(9.7,9.2)(10.3,9.8)
 \psline[linestyle=dotted](7,3.5)(9,4.5)
 \psline[linestyle=dotted](7,7.5)(9,8.5)
 \psline[linestyle=dotted](9,2.5)(7,8.5)
\end{pspicture}

\textit{Cojacobi}: The proof of Cojacobi for $s_2$ uses several
versions of diagram IV. Again the direction on the chords is
related to the order of arguments and subsequent orientations.
Each column represents $(s_2 \otimes 1)\cdot s_2$ or its cyclic
permutation. The numbers indicate the output arguments. The
particular 4 element subset $F$ pictured with its indicated
coincidences could contribute in various ways to the chain
operation.  A chord could have either direction or be used as
first chord or as second chord in the composed operation. Of the 8
possibilities in each case the 4 pictured are the only ones that
occur. The cancellation is indicated.

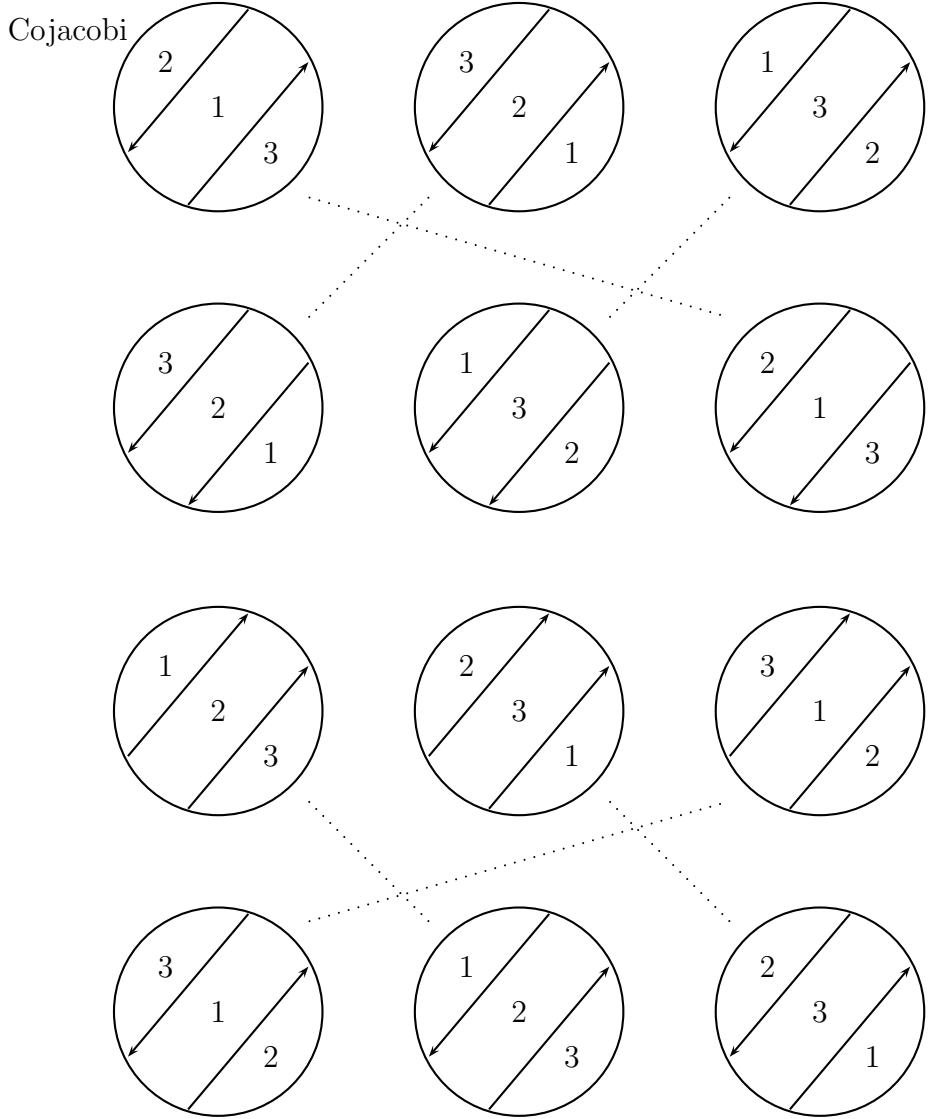
\begin{pspicture}(0,0)(12,8)
%\psgrid(0,0)(12,8)
 \pscircle(2,2){1.4} \pscircle(6,2){1.4} \pscircle(10,2){1.4}
 \pscircle(2,6){1.4} \pscircle(6,6){1.4} \pscircle(10,6){1.4}
 \psline{<-}(1.6,0.7)(3.2,2.6) \psline{<-}(0.8,1.4)(2.4,3.3)
 \rput(0,7){Cojacobi}
 \rput(1.3,2.6){$3$} \rput(2,2){$2$} \rput(2.7,1.4){$1$}
 \psline{<-}(5.6,0.7)(7.2,2.6) \psline{<-}(4.8,1.4)(6.4,3.3)
 \rput(5.3,2.6){$1$} \rput(6,2){$3$} \rput(6.7,1.4){$2$}
 \psline{<-}(9.6,0.7)(11.2,2.6) \psline{<-}(8.8,1.4)(10.4,3.3)
 \rput(9.3,2.6){$2$} \rput(10,2){$1$} \rput(10.7,1.4){$3$}
 \psline{->}(1.6,4.7)(3.2,6.6) \psline{<-}(0.8,5.4)(2.4,7.3)
 \rput(1.3,6.6){$2$} \rput(2,6){$1$} \rput(2.7,5.4){$3$}
 \psline{->}(5.6,4.7)(7.2,6.6) \psline{<-}(4.8,5.4)(6.4,7.3)
 \rput(5.3,6.6){$3$} \rput(6,6){$2$} \rput(6.7,5.4){$1$}
 \psline{->}(9.6,4.7)(11.2,6.6) \psline{<-}(8.8,5.4)(10.4,7.3)
 \rput(9.3,6.6){$1$} \rput(10,6){$3$} \rput(10.7,5.4){$2$}
 \psline[linestyle=dotted](3.2,3.2)(4.8,4.8)
 \psline[linestyle=dotted](7.2,3.2)(8.8,4.8)
 \psline[linestyle=dotted](3.2,4.8)(8.8,3.2)
\end{pspicture}

\begin{pspicture}(0,0)(12,8)
%\psgrid(0,0)(12,8)

\pscircle(2,2){1.4} \pscircle(6,2){1.4} \pscircle(10,2){1.4}
 \pscircle(2,6){1.4} \pscircle(6,6){1.4} \pscircle(10,6){1.4}
 \psline{->}(1.6,0.7)(3.2,2.6) \psline{<-}(0.8,1.4)(2.4,3.3)
 \rput(1.3,2.6){$3$} \rput(2,2){$1$} \rput(2.7,1.4){$2$}
 \psline{->}(5.6,0.7)(7.2,2.6) \psline{<-}(4.8,1.4)(6.4,3.3)
 \rput(5.3,2.6){$1$} \rput(6,2){$2$} \rput(6.7,1.4){$3$}
 \psline{->}(9.6,0.7)(11.2,2.6) \psline{<-}(8.8,1.4)(10.4,3.3)
 \rput(9.3,2.6){$2$} \rput(10,2){$3$} \rput(10.7,1.4){$1$}
 \psline{->}(1.6,4.7)(3.2,6.6) \psline{->}(0.8,5.4)(2.4,7.3)
 \rput(1.3,6.6){$1$} \rput(2,6){$2$} \rput(2.7,5.4){$3$}
 \psline{->}(5.6,4.7)(7.2,6.6) \psline{->}(4.8,5.4)(6.4,7.3)
 \rput(5.3,6.6){$2$} \rput(6,6){$3$} \rput(6.7,5.4){$1$}
 \psline{->}(9.6,4.7)(11.2,6.6) \psline{->}(8.8,5.4)(10.4,7.3)
 \rput(9.3,6.6){$3$} \rput(10,6){$1$} \rput(10.7,5.4){$2$}
 \psline[linestyle=dotted](3.2,4.8)(4.8,3.2)
 \psline[linestyle=dotted](7.2,4.8)(8.8,3.2)
 \psline[linestyle=dotted](3.2,3.2)(8.8,4.8)
\end{pspicture}

\textit{Drinfeld compatibility}: The proof of Drinfeld
compatibility between $c_2$ and $s_2$ uses diagram V and diagram
VI. The figure represents $s_2[1,2]$. The directions on chords
indicate ordering of input variables for $s_2[1,2]$ and subsequent
orientations. The dot on the chord indicates it is used first.
Cancellation is indicated. The last four terms represent the right
hand side of the compatibility equation,
$s_2([1,2])=[s_2(1),2]+[1,s_2(2)]$:

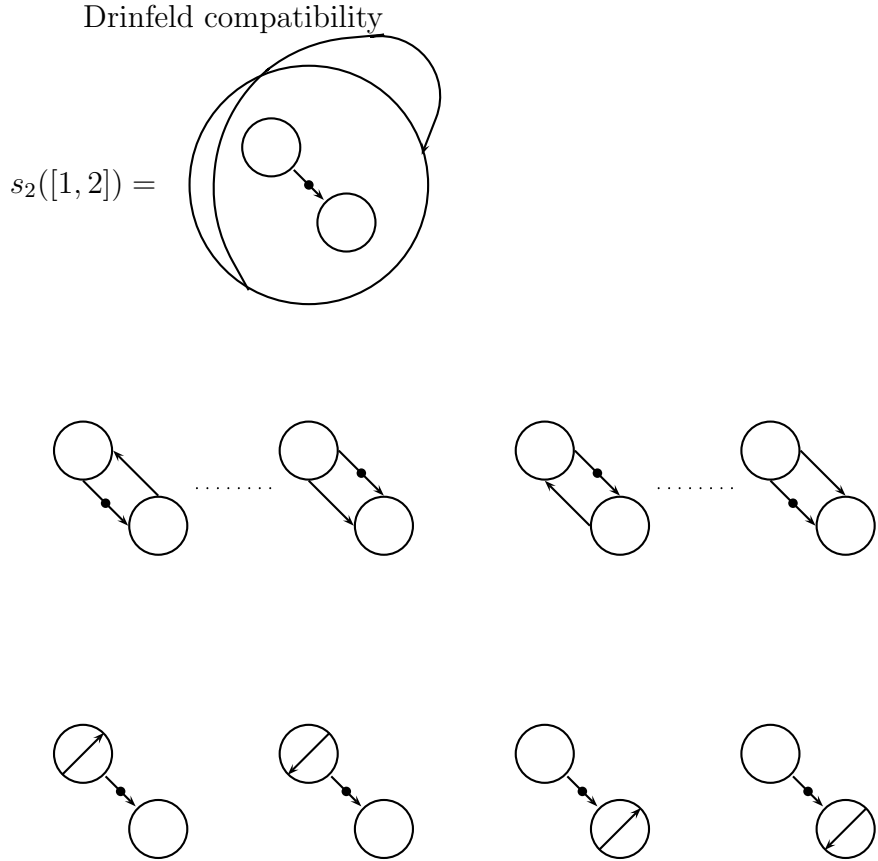
\begin{pspicture}(-2,0)(4,4.5)
 %\psgrid(-2,0)(4,4.5)
 \pscircle(2,2){1.6}
 \pscircle(1.5,2.5){0.4}
 \pscircle(2.5,1.5){0.4}
 \psline{->}(1.8,2.2)(2.2,1.8) \psdot(2,2)
 \psline[linearc=2](1.2,0.6)(0.7,1.5)(0.8,3)(1.9,3.9)(3,4)
 \psline[linearc=0.8]{->}(3,4)(4.2,4.2)(3.5,2.4)
 \rput(-1,2){$s_2 ([1,2])=$}
 \rput(1,4.2){Drinfeld compatibility}
\end{pspicture}

\begin{pspicture}(0.5,0)(6.5,4)
 %\psgrid(0.5,0)(6.5,4)
 \pscircle(1.5,2.5){0.4}
 \pscircle(2.5,1.5){0.4}
 \pscircle(4.5,2.5){0.4}
 \pscircle(5.5,1.5){0.4}
 \psline[linestyle=dotted](3,2)(4,2)
 \psline{->}(1.5,2.1)(2.1,1.5)
 \psline{->}(2.5,1.9)(1.9,2.5) \psdot(1.8,1.8)
 \psline{->}(4.5,2.1)(5.1,1.5)
 \psline{<-}(5.5,1.9)(4.9,2.5) \psdot(5.2,2.2)
\end{pspicture}
\begin{pspicture}(0.5,0)(6.5,4)
 %\psgrid(0.5,0)(6.5,4)
 \pscircle(1.5,2.5){0.4}
 \pscircle(2.5,1.5){0.4}
 \pscircle(4.5,2.5){0.4}
 \pscircle(5.5,1.5){0.4}
 \psline[linestyle=dotted](3,2)(4,2)
 \psline{<-}(1.5,2.1)(2.1,1.5)
 \psline{<-}(2.5,1.9)(1.9,2.5) \psdot(2.2,2.2)
 \psline{->}(4.5,2.1)(5.1,1.5)
 \psline{<-}(5.5,1.9)(4.9,2.5) \psdot(4.8,1.8)
\end{pspicture}

\begin{pspicture}(0.5,0)(6.5,4)
% \psgrid(0.5,0)(6.5,4)
 \pscircle(1.5,2.5){0.4}
 \pscircle(2.5,1.5){0.4}
 \psline{->}(1.8,2.2)(2.2,1.8) \psdot(2,2)
 \psline{->}(1.22,2.22)(1.78,2.78)
 \pscircle(4.5,2.5){0.4}
 \pscircle(5.5,1.5){0.4}
 \psline{->}(4.8,2.2)(5.2,1.8) \psdot(5,2)
 \psline{<-}(4.22,2.22)(4.78,2.78)
\end{pspicture}
\begin{pspicture}(0.5,0)(6.5,4)
%\psgrid(0.5,0)(6.5,4)
 \pscircle(1.5,2.5){0.4}
 \pscircle(2.5,1.5){0.4}
 \psline{->}(1.8,2.2)(2.2,1.8) \psdot(2,2)
 \psline{->}(2.22,1.22)(2.78,1.78)
 \pscircle(4.5,2.5){0.4}
 \pscircle(5.5,1.5){0.4}
 \psline{->}(4.8,2.2)(5.2,1.8) \psdot(5,2)
 \psline{<-}(5.22,1.22)(5.78,1.78)
\end{pspicture}

 \textit{Identities for $c_n$ and $s_n$}: Using
computations of Getzler [G] one can show the $c_n$ taken together
satisfy a generalized Jacobi identity.  Some but not all of these
were shown at the chain level in [CS]. The rest follow as in [G]
from the BV homology structure of [CS]. The Getzler identities
(defining what he calls a "gravity algebra") are Koszul dual [G]
to the associative (or commutative) identities in the definition
of a "Frobenius manifold" [M] describing the algebraic structure
of genus zero Gromov Witten invariants of a closed symplectic
manifold.

In the latter case there is a compatible non degenerate inner
product- the Poincare duality of a closed symplectic manifold in
its Floer homology which is equal to its ordinary homology.  One
can always form the "positive boundary" version of an algebraic
structure by trading in multiplications and inner product for
mulitplications and co multiplications satisfying the induced
identities.  (This idea we learned from David Kahzdan and the
terminology "positive boundary" was suggested by Ralph Cohen.)

Do this for "Frobenius manifold", form the positive boundary
version and then apply Koszul duality [G], [Gan], [M, p.87].  One
obtains an algebraic structure which combines gravity algebra and
gravity co algebra with Drinfeld type compatibilities.  It
contains the notion of Lie bialgebra [Gan].  At this point one
knows the $c_n$ satisfy the gravity algebra identities
(generalized Jacobi [G]), the $s_n$ by arrow reversal duality
satisfy the gravity coalgebra identities (generalized Cojacobi)
and $c_2$ and $s_2$ satisfy Drinfeld compatibility (this paper).

This is the evidence for the conjecture and problem mentioned
above.

 \textit{Involutive property of the Lie bialgebra}: We may consider
the operation $e=c_2 \circ s_2$:
$\mathbb{H}\rightarrow\mathbb{H}$. If we think of $s_2$ and $c_2$
as associated to pairs of pants pointed in opposite directions
then $e$ is associated to a torus with one input circle and one
output circle obtained by glueing these two pairs of pants. The
diagram for $e$ is one circle with (ordered) two chords whose
endpoints are linked (see diagram VII). The (relative)
compactification of this Diagram VII uses all of the
considerations above- the different parts colliding of Diagrams
III and IV, the multiple points of Diagram V, and the FM stratum
of Diagram VI.  Examining the chain operation shows that there is
a complete cancellation because interchanging the ordering of
chords is orientation reversing for Diagram VII. Thus the genus
one operator $e$ is always zero in the Lie bialgebra.

Now $e$ is the infinitesimal analogue for Lie bialgebras of the
square of the canonical antiautomorphism of a Hopf algebra.

When this square is the identity one says the Hopf algebra is
involutive so we say in analogy that since $e$ is zero the
\textit{Lie bialgebra of the Theorem is involutive.} Note this is
a genus one relation.

\textit{Higher genus operators}: Forming the (relative)
compactification of a general chord diagram $D$ requires the
addition of many strata.  When different parts of $F$ collide at
one point the normal bundles fit together as indicated above. When
the same part self collides at one point we introduce a
multiplicity at that point.  When several points collide at
different rates we can reduce to the previous cases.  When several
points collide at commensurable rates we may have to introduce FM
strata as in paragraph Diagram VI.  The constraint bundle issue
must be solved by adding relations among derivatives at the
coincident points.  These relations are associated to cycles in
the subgraphs of $\Gamma (D)$ which are collapsing. For example
the multiple point of Diagram V (and other examples) corresponds
to collapsing a loop of $\Gamma(D)$.  The FM stratum of Diagram VI
corresponds to collapsing a cycle in $\Gamma (D)$ made out of two
edges.  More generally when intervals between points of $F$
collapse this determines a collapsing subgraph of $\Gamma(D)$.
Cycles on this graph give relations among derivatives.  These are
required to define the constraint normal bundle.

The theory begins to take on the structure of the collapsing
graphs in the renormalization theory of Feynman diagrams appearing
in the work of Kreimer et al [K].  This will be discussed
elsewhere [S].

\newpage
\textit{References}:

[cs] M. Chas and D. Sullivan "String Topology" GT/ 9911159.
Accepted by Annals of Mathematics.

[V] A. Voronov "Notes on Universal Algebra". See author's web
page.  To appear in proceedings of Stony Brook conference on
"Graphs and Patterns in Mathematics and Physics", June 2001.

[G] E.Getzler "Operads and moduli spaces of genus zero Riemann
surfaces".  In:The Moduli Spaces of Curves ed. by R.Dijkgraaf, C.
Faber, G. van der Geer Progress in Math, vol.129 Birkhauser 1995,
199-230.

[Gan] Wee Liang Gan "Koszul duality for dioperads" preprint
University of Chicago 2002. QA/0201074

[M] Yuri I. Manin "Frobenius Manifolds, Quantum cohomology and
Moduli spaces". AMS Colloquium Publications Vol.47.

[D] Pierre Deligne "Quantum Fields and Strings: A course for
mathematicians" vol. 1 chapter 1 section 1.2 Editors Deligne et al
AMS IAS

[FM] W. Fulton, R. MacPherson "A compactification of configuration
spaces" Annals of Math, 139 (1994) 183-225

[K] D. Kreimer "Structures in Feynman graphs- Hopf algebras and
symmetries"  See author's web page.  To appear in proceedings of
Stony Brook conference on "Graphs and Patterns in Mathematics and
Physics", June 2001.

[S] D. Sullivan, work in progress "Closed String Operators and
Feynman graphs", possible title.

[T] V. Turaev "Skein quantization of Poisson algebra of loops on
surfaces" Ann. Sci. Ecole Norm. Sup. (1) 24 (1991), no.6, 635-704.

[MC] Moira Chas "Combinatorial Lie Bialgebras of curves on
surfaces" GT/0105178, to appear in Topology.
\end{document}